# FRACTIONS IN THE *SUÀN SHÙ SHŪ*
# (CHINA, BEGINNING OF THE 2ND CENTURY BCE)

# 西漢出土文獻《算數書》分數表達

**Rémi Anicotte 安立明**
Associate member of CRLAO, Paris
[法] 東亞語言研究所准成員

**ABSTRACT**：The *Suàn Shù Shū* contains 301 instances of regular expressions for fractions. They can be "mono-dimensional" (formed with one integer name only) for unit fractions, "bidimensional" (with two integer names) for both unit and non-unit fractions, or lexicalized only for 1/3, 1/2 and 2/3. The present paper gives a complete description of the diversity of these forms. Bidimensional expressions are predicative phrases: the name *n fēn* of a unit fraction 1/*n* acts as subject and the numerator's name as predicate; according to the syntactic context, the morpheme *zhī* can be used as an optional marker of this predicative relation.

**KEYWORDS**：Chinese historical syntax, Fractions, Numerical expressions, Separable semantic units, Measure words, *Suàn Shù Shū*.

摘要：西漢出土文獻《算數書》中表達分數的短語有 301 例。本文全面描述該書中的那些短語。在文中，我把它們定義為"一維短語（即包括唯一一個整數名）"和"二維短語（即包括兩個整數名）"，前者表達單位分數，後者表達單位和非單位分數，除此之外還有三個詞彙化的短語專門用於表達 1/3、1/2 和 2/3。在二維短語中，單位分數 1/*n* 作為短語的主語表達分母，而表達分子的那個數名作為短語的謂語，根據語境，"之"字可作為短語標誌放在主謂語之間。

主題詞：中文語法歷史，分數，表数短語，離合詞組，量詞，《算數書》。

This is a preprint version of a paper published in *The Journal of Chinese Linguistics*, vol.45 no.1 (January 2017): 20-67. The author teaches at Beijing International French School and can be contacted by email at *remi_anicotte@yahoo.com*. He wishes to thank Karine Chemla (SPHERE), Redouane Djamouri (CRLAO), Christoph Harbsmeier (University of Oslo), Sylviane Schwer (Paris 13 University), Xiao Can (Hunan University), Xu Dan (INALCO, IUF) and the reviewers of the Journal of Chinese Linguistic for their invaluable encouragement and help with the analysis of the corpus.

Abbreviations DECL: declarative; MW: measure word; NUM: numeral; OBJ: object of transitive verb; 3OBJ: third person object pronoun; {*n*} (with a number *n* written in Arabic numerals): the mono-morphemic expression of the number *n* in a given language; *A*(*B*) and *A*(*BC*): the character *A* is a rendition of the original character encountered in the Chinese corpus, the character *B* or the sequence *BC* in parenthesis is a modern form of what is understood for *A*. For example: 有(又), 廿(二十), 卅(三十), 卌(四十), 丰(七十), 泰(大).

## 1. FRACTIONS IN THE *SUÀN SHÙ SHŪ*

The *Suàn Shù Shū*[1] is a mathematical text written on 190 bamboo strips, consisting of approximately 7,000 characters. The manuscript was excavated during the winter of 1983–84 from a Han Dynasty tomb in Zhangjiashan where a calendar for the year 186 BCE was found, and so the tomb is thought to have been closed that very year: the book was probably written in the beginning of the 2$^{nd}$ century BCE. Peng Hao (2001: 4–6) states that the production and taxation standards mentioned in some passages prove that they were copied from originals written in the kingdom of Qin before the unification of China in 221 BCE, while other sections could only have been composed during the reign of the Western Han Dynasty which began in 206 BCE.

The text contains 301 occurrences (not all different) of regular expressions for fractions. This is much more than in Qin-Han manuscripts found in the same tomb[2], or elsewhere but not specialized in mathematics. More Qin-Han texts of mathematics were discovered after the *Suàn Shù Shū*; they[3] shows no discrepancy concerning the expression of integers and fractions which are all *proper fractions*, i.e. smaller than 1, the denominator being larger than the numerator[4]. They are all written in the Chinese language; I use Arabic numerals in translations because they are more readable than numbers written in English, but there is no symbolic numerical notation in the original.

Expressions for fractions can be special lexicalized items, but clearly such a scheme is viable only for a few specific fractions; in Qin-Han texts this was limited to 1/2, 1/3 and 2/3. On the other hand a *generic* linguistic pattern capable of communicating the fraction of any two integers has to account for both the numerator and the denominator, thus producing numerical expressions which I call *bidimensional*. Unit fractions 1/$n$ are a special case in the corpus since approximately 64% (exactly 83 out of 129) of them only state their denominator $n$ and end up *monodimensional*.

English *two-fifths* can be inserted before nouns or measure words as the head of a noun phrase using the preposition *of*; for example *two-fifths of the*

---

[1] Authorised editions of the *Suàn Shù Shū* were successively published in [Wenwu 2000] and [Wenwu 2001]. An edition with commentary was published by Peng Hao (2001). Japanese translations were produced by Jochi S. (2001) and [Ōkawa et al. 2006]; the latter includes a translation into Contemporary Chinese by Ma Biao. Another Contemporary Chinese rendition was done by Hu Yitao (2006) who worked under the supervision of Zhang Xiancheng (Southwest University, Chongqing). Two independent English translations were successively produced by Cullen (2004) and Dauben (2008).

[2] Yang Lingrong (2008: 14-20) counts 328 fractions in the *Suàn Shù Shū* (including atypical expressions of fractions and expressions of proportions), but 58 fractions only in all the other corpora of Zhangjiashan.

[3] I checked *Shù* 數 (it belonged to the Yuèlù Academy, the text is now available in Xiao Can 2010; Xiao Can briefed me by email about integers and fractions in *Shù* before the formal publication), and *Suàn Shù* 算術 (which was excavated in Shuihudi in 2006; at the time of my research, only short excerpts were available in Xiong Beisheng et al. (2008), a joint publication by the Hubei Provincial Institute of Cultural Relics and Archaelogy (Húběi Shěng Wénwù Kǎogǔ Yánjiūsuǒ) and the Yunmen Museum (Yúnmèng Xiàn Bówùguǎn) published in 2008, and in Chemla & Ma Biao (2011)).

[4] There is a wide range of values in the *Suàn Shù Shū* that depend on the calculations they are involved in. For example, 1/50 is on strip 4, 12/18 is on strip 55, 47/98 is on strip 92, and 162/2016 is on strip 20.

*population* or *two-fifths of a liter*. Contemporary Chinese expresses 2/5 as *wǔ fēn zhī èr*, i.e. {5} *fēn zhī* {2}, with the denominator's name given first, the compound *fēn zhī* coming next and the numerator's name in last position. The compound {5} *fēn zhī* {2} can in turn be inserted directly before a measure word or a noun to form the sequences "Fraction Name + MW" or "Fraction Name + Noun" respectively; placing the fraction name in the head position of a noun phrase "Noun (+ *de* 的) + Fraction Name" is also possible. Nothing can be introduced between the components of *two-fifths* or {5} *fēn zhī* {2}, therefore such compounds are *inseparable semantic units* or *inseparable fraction names*. But in Qin-Han manuscripts, the only inseparable fraction names were on one hand special lexicalized expressions of 1/3, 1/2 and 2/3, and on the other hand the monodimensional expressions of unit fractions built according to the pattern "Denominator + *fēn*"; measure words were inserted after these expressions. Bidimensional expressions of fractions were built as predicative phrases with the name "Denominator's Name + *fēn*" of a unit fraction acting as subject and with the numerator's name acting as predicate. When a measure word was involved it was inserted right after "Denominator + *fēn*". The morpheme *zhī* could be optionally added before the numerator's name as a marker of the predicative relation.

The present paper provides an exhaustive survey of the diversity of all the expressions for fractions in the *Suàn Shù Shū*[5].

## 2. MEASURE WORDS, INTEGERS AND MIXED NUMBERS

Measure words can be found after the "Denominator + *fēn*" compounds. Measure words in the *Suàn Shù Shū* are mostly units of measurement[6]; they fit into the construction "NUM + MW" where the numeral can be the name of an integer, as well as lexicalized fraction names and monodimensional unit fraction expressions.

---

[5] Guo Shuchun (2002) and Yang Lingrong (2008) already presented the various patterns which can be encountered, but they failed to quantify their distribution and to relate them to their insertion contexts.

[6] Length units are *cùn* 寸, *chǐ* 尺, *bù* 步, *zhàng* 丈 and *lǐ* 里: 1 *lǐ* = 180 *zhàng*, 1 *zhàng* = 10 *chǐ*, 1 *bù* = 6 *chǐ*, 1 *chǐ* = 10 *cùn*. There is also a specific unit *wéi* 韋(圍) used only for circumferences according to Peng Hao (2001): 1 *wéi* = 3 *chǐ* (i.e. ≈ π *chǐ*, the circumference of a circle with diameter 1 *chǐ*). Surface area units can be *mǔ* 畝 and *qǐng* 頃 (1 *qǐng* = 100 *mǔ*), or are derived from length units and contextually understood as referring to surfaces even though there is no special indication equivalent to English "square". Capacity units are *shēng* 升, *dǒu* 斗 and *shí* 石 (1 *shí* = 10 *dǒu*, 1 *dǒu* = 10 *shēng*). The only volume unit in the *Suàn Shù Shū* is *chǐ* 尺; it is derived from the length unit *chǐ* 尺 and occurs without any special indication equivalent to "cubic". Weight units are *zhū* 朱(銖), *liǎng* 兩, *jīn* 斤, *jūn* 鈞 and *shí* 石: 1 *shí* = 4 *jūn*, 1 *jūn* = 30 *jīn*, 1 *jīn* = 16 *liǎng*, 1 *liǎng* = 24 *zhū*. Note that *shí* 石 can designate both a capacity unit and a weight unit, and that the *Xiàndài hànyǔ guīfàn cídiǎn* 现代汉语规范词典 [Dictionary of contemporary Chinese] (Beijing, 2010: 262) gives the pronunciation *dàn* for this character used as a measurement unit, but recommends *shí* when reading ancient texts.

Mass nouns *sù* 粟 [unhusked millet][7], *mǐ* 米 [husked millet], *bài* 粺 [milled millet][8], *shuǐ* 水 [water], *qī* 柒(漆) [lacquer], *jīn* 金 [gold], *guǎng* 廣 [width], *zòng* 縱 [length], etc. occur in "Noun + NUM + MW" sequences, whereas *rén* 人 [person] and nouns for countable items like *lútáng* 盧唐 [bamboo tube], *jiǎn* 簡 [bamboo strip], *suàn* 筭(算) [string of coins][9], etc. all fit into the pattern "NUM + Noun" in the same manner as measure words[10]. The word *qián* 錢 is used as a currency unit; it occurs either in "NUM + MW" or in "Noun + NUM". The numeral 1 is not always stated before a measure word; this is marked with ∅ right before the measure word *chǐ* in (9), (42), and (43), and *wéi* in (69).

Fraction expressions are built from the names of the numerator and denominator, which are integers. Chinese names for integers belong to a decimal numeration. The digits are *yī* 一 {1}[11], *èr* 二 {2}[12], *sān* 三 {3}, *sì* 四 {4}, *wǔ* 五 {5}, *liù* 六 {6}, *qī* 七 {7}, *bā* 八 {8}, *jiǔ* 九 {9}. In the *Suàn Shù Shū*, the series of pivots is limited to *shí* 十 {10}, *bǎi* 百 {$10^2$}, *qiān* 千 {$10^3$} and *wàn* 萬 {$10^4$}; the largest number being 10,000,000 expressed as *qiān wàn*, i.e. {$10^3$}{$10^4$}, on strip 11.

I use the notation {*number*} with a number written in Arabic numerals between braces to represent the numerical morpheme which expresses the bracketed number in a particular language. For instance the same notation {10} can represent the word *ten* in English and the morpheme *shí* in Chinese. The notation {$10^4$} represents *wàn* in Chinese, but would not occur for English *ten thousand*, which is represented as {10}{$10^3$} and stands for the succession of the mono-morphemic items *ten* {10} and *thousand* {$10^3$}.

In the *Suàn Shù Shū*, the morpheme {1} is used before all pivots in a number name but the highest one. This is visualized with a shaded *yī* 一 {1} in examples (1), (2), (4), (6), and a shaded ∅ in (3), (5), (6), (7). The sequences {*digit*}{*pivot*} and {*smaller pivot*}{*larger pivot*} express products, they are concatenated directly. These concatenations express sums.

The conjunction *yòu* was used in the Western Zhou inscriptions on bronze vessels to join the tens and units places, and sometimes also the hundreds and tens, but in the *Suàn Shù Shū* it occurs only in expressions of mixed numbers to link an integer and a fraction as in (7) and it is free pattern[13]. The term *líng* 零 is not encountered anywhere in the corpus and was not used in integer names before the 12th century CE[14].

---

[7] I refer to Dauben (2008: 169-170) for the English names of crops and their byproducts in the *Suàn Shù Shū*.

[8] In the *Suàn Shù Shū* as in the *Nine Chapters* [*Jiǔ Zhāng Suàn Shù* 九章算術], *bài* 粺 refers to milled millet, not to a particular variety of millet.

[9] The word *suàn* 筭(算) refers either to a string of coins (Cullen 2004: 29) or to a Han dynasty unit of taxation (Chemla & Guo Shuchun 2004: 989).

[10] In Chinese the distinction between mass and countable nouns is semantic, not grammatical.

[11] According to the rules of *pīnyīn* transcription, the digit {1} is always romanized as *yī* with a first tone mark regardless of the actual tone in Contemporary Chinese. This depends on the tone of the following syllable; for example {1}{$10^4$} is actually pronounced *yí wàn*, but this is only noted in publications concerned with pronunciation.

[12] The Contemporary Chinese variant *liǎng* 兩 was hardly used in exact number names before the 20th century CE and the digit *èr* 二 is the only numeral for 2 in the *Suàn Shù Shū*.

Additionally, to write tens, the copyists of the *Suàn Shù Shū* could use ligatures instead of the corresponding two-character forms. The numbers 50, 60, 80 and 90 were either written with two separate characters, or as ligatures of the characters for 5, 6, 8 or 9 written in a reduced size above the character for 10; I always transcribe these with the two-character forms *wǔ shí* 五十, *liù shí* 六十, *bā shí* 八十 and *jiǔ shí* 九十 respectively. But 20, 30, 40 and 70 were always written 廿, 卅, 卌 and 〧[15] respectively. Conforming to the current scholarly usage, I transcribe them with their actual written forms followed by the two-character forms in parentheses: *èr shí* 廿(二十), *sān shí* 卅(三十), *sì shí* 卌(四十) and *qī shí* 〧(七十) respectively; the *pīnyīn* romanization and glosses are those of the disyllabic compounds[16]; see for example the expressions with 20 in (4) and (6).

(1)  
in *Suàn Shù Shū* strip 172  
二百　　一十  
*èr bǎi　　yī shí*  
$\{2\}\{10^2\}　\{1\}\{10\}$  
'210'

(2)  
in *Suàn Shù Shū* strip 20  
二千　　一十　　六  
*èr qiān　　yī shí　　liù*  
$\{2\}\{10^3\}　\{1\}\{10\}　\{6\}$  
'2016'

(3)  
in *Suàn Shù Shū* strip 76  
錢　　∅百　　五十  
*qián　　∅bǎi　　wǔ shí*  
*qián*　$∅\{10^2\}　\{5\}\{10\}$  
'150 *qián*'

(4)  
in *Suàn Shù Shū* strip 176  
七千　　一百　　廿(二十)　　九  
*qī qiān　　yī bǎi　　èr shí　　jiǔ*  
$\{7\}\{10^3\}　\{1\}\{10^2\}　\{2\}\{10\}　\{9\}$  
'7129'

---

[13] When a measure word separates the integer and the fraction, the item *yòu* can be present as in (7) or absent as in (18). When there is no measure word, the item *yòu* is sometimes not used as in (27) or used as in (119).

[14] Readers can find more details on Chinese integer names in Anicotte (2015 a).

[15] The ligature 〧 for 70 is also found in the Qin-Han manuscripts *Shù* and *Suàn Shù* mentioned above, but not in dictionaries of Middle Chinese.

[16] Modern dictionaries give the pronunciations *niàn* for 廿, *sà* for 卅, and *xì* for 卌, but other ligatures for tens are never mentioned; this justifies the present usage of glossing *all* of these ligatures as dissyllabic compounds. I presume nothing about the actual pronunciation of tens during the Qin-Han period; this matter is beyond the scope of the present paper.

| | | | | | |
|---|---|---|---|---|---|
| (5) | ∅千 | 八十 | 九 | | |
| in *Suàn Shù Shū* | ∅*qiān* | *bā shí* | *jiǔ* | | |
| strip 172 | ∅{$10^3$} | {8}{10} | {9} | | |
| | '1089' | | | | |
| (6) | ∅萬 | 一千 | 五百 | 廿(二十) | 銖 |
| in *Suàn Shù Shū* | ∅*wàn* | *yī qiān* | *wǔ bǎi* | *èr shí* | *zhū* |
| strip 47 | ∅{$10^4$} | {1}{$10^3$} | {5}{$10^2$} | {2}{10} | *zhū* |
| | '11520 *zhū*' | | | | |
| (7) | ∅十六 | 尺 | 有(又) | 十八分 | 尺 |
| in *Suàn Shù Shū* | ∅*shí liù* | *chǐ* | *yòu* | *shí bā fēn* | *chǐ* |
| strip 55 | ∅{10}{6} | *chǐ* | and | {10}{8} *fēn* | *chǐ* |
| | '16 *chǐ* and 12/18 *chǐ*'[17] | | | | |

∅十二
∅*shí èr*
∅{10}{2}

The mixed number 16 12/18 in (7) is inserted with the measure word *chǐ*. The integer and the fraction are dealt with as two independent quantification phrases which are concatenated with *yòu* here or juxtaposed in other instances. The measure word *chǐ* occurs first with the integer and is inserted again within the expression for the fraction 12/18. This is a regular pattern in the *Suàn Shù Shū*[18], even the nouns *jiǎn* 簡 [bamboo strip] and *suàn* 筭(算) [string of coins], in (55) and (57) respectively, appear first after the name of the integer and then again with the proper fraction after *fēn*. There are only a few exceptions: 3 can be seen in (28)–(30) with a measure word stated only after the integer, and there is also (76) with the noun *lútáng* 盧唐 [bamboo tube] placed in the measure word position in the integer expression but not repeated after *fēn*. On the other hand, (31) is a unique case of over-repetition, with the measure word *bù* also found after the denominator's name and therefore appearing three times.

### 3. GENERIC UNIT FRACTION EXPRESSIONS
### 3–1. Monodimensional unit fraction expressions

In the *Suàn Shù Shū* as in other Qin-Han texts of mathematics, the canonical names for unit fractions 1/*n* were *n fēn*. They stated only the denominator *n* but not the numerator 1 and were therefore monodimensional numerical expressions. Nothing could be inserted between the two constituents *n* and *fēn*; that is, the

---

[17] English would simply write 16 12/18 *chǐ* and say "sixteen and twelve-eighteenths of a *chǐ*". My translations are meant to highlight the repetition of the measure words and the presence of linking terms in the original.

[18] This repetition shows that the integer and the fraction were dealt with independently; this has nothing to do with "echo" constructions which occur in *one* quantification phrase "Noun + Num + Noun" (e.g. *qiāng bǎi qiāng* 羌百羌 [Qiang hundred Qiang] i.e. 'one hundred Qiang people' on bone inscription H32042).

*n fēn* compounds were inseparable semantic units. Measure words followed *n fēn* just as they followed integers. In other words, these *n fēn* unit fraction designations were numerals in their own right.

The item *fēn* in *n fēn* can be considered syntactically neutralized and fossilized through its involvement in a word-formation process. It was used on occasion in the corpus as a noun for *part* or *fraction,* in which cases it could be read *fèn* with a falling tone regardless of the modern reading *fēn*. It was also used as a verb meaning *to share* (which should be read as *fēn*); but considering *fēn* in the compounds *n fēn* synchronically, it seems pointless to try and interpret it as a noun or as a verb. From a semantic point of view, the question would be whether *fēn* referred more to the action of partitioning or more to the result (the parts) of this action. This might be interesting in an attempt to reconstruct the *emergence* of the expressions *n fēn*, but again this seems irrelevant synchronically.

In the *Suàn Shù Shū*, there are 83 instances (not all different) of monodimensional expressions of 1/$n$ unit fractions. Among them, 76 do not involve a measure word, while 8 instances do. The 76 instances without measure words are:

> → *sān fēn* for 1/3 (16 instances, on strips 3, 5, 6, 24, 27, 139, 168, 169, 170, 171, 172, 174, 176 and 179);
> → *sì fēn* for 1/4 (15 instances, on strips 3, 5, 6, 9, 14, 27, 169, 171, 172, 174, 176 and 179);
> → *wǔ fēn* for 1/5 (14 instances, on strips 5, 6, 8, 9, 14, 170, 171, 172, 174, 176 and 179);
> → *liù fēn* for 1/6 (10 instances, on strips 9, 10, 171, 172, 174, 176 and 179);
> → *qī fēn* for 1/7 (9 instances, on strips 5, 9, 10, 172, 174, 176 and 179);
> → *bā fēn* for 1/8 (5 instances, on strips 10, 174, 176 and 179);
> → *jiǔ fēn* for 1/9 (2 instances, on strips 176 and 177);
> → *shí fēn* for 1/10 (2 instances, on strips 179 and 180);
> → *bǎi fēn* for 1/100 (2 instances, on strips 14 and 16).

The 8 instances of expressions for unit fractions followed by a measure word are:

> → *sān fēn cùn* for 1/3 *cùn* on strip 2;
> → *bā fēn cùn* for 1/8 *cùn* on strip 2;
> → *sì fēn cùn* for 1/4 *cùn* on strip 4;
> → *wǔ fēn cùn* for 1/5 *cùn* on strip 4;
> → *liù fēn cùn* for 1/6 *cùn* on strip 4;
> → *liù shí fēn chǐ* for 1/60 *chǐ* on strip 4;
> → *jiǔ fēn zhū* for 1/9 *zhū* on strip 29.

**3–2. Bidimensional expressions of unit fractions**

The numerator's name {1} is not compulsory and is usually omitted when the fraction comes as a factor in a multiplication: check ∅ right after *fēn* in (8) with no measure word, and after the measure word *cùn* in (9).

(8)　　　　　　四分⌀　　乘　　　四分⌀　　十六分　　　一
in *Suàn Shù Shū*　*sì fēn* ⌀　*chéng*　*sì fēn* ⌀　*shí liù fēn*　*yī*
strip 9　　　　　{4} *fēn* ⌀　multiply　{4} *fēn* ⌀　{10}{6} *fēn*　{1}
　　　　　　　'1/4 times 1/4 are 1/16'

(9)　　　　　　五分　　寸⌀　　乘　　　　⌀尺
in *Suàn Shù Shū*　*wǔ fēn*　*cùn* ⌀　*chéng*　　⌀ *chǐ*
strip 4　　　　　{5} *fēn*　*cùn* ⌀　multiply　⌀ *chǐ*
　　　　　　　'1/5 *cùn* times [one] *chǐ* are

　　　　　　　五十分　　　尺　　一　　也
　　　　　　　*wǔ shí fēn*　*chǐ*　*yī*　*yě*
　　　　　　　{5}{10} *fēn*　*chǐ*　{1}　DECL
　　　　　　　1/50 [square] *chǐ*'[19]

But {1} is not omitted (yielding bidimensional numerical expressions) when stating the result of multiplications involving *n fēn* unit fraction names or lexicalized names for 1/2 and 1/3. Examples can be seen in (8) and (9) above for the former situation, and in (117) and (122) from Sect. 5 for the latter one. We can analyze the combination of *n fēn* and {1} as a predicative clause with the monodimensional name of the unit fraction acting as the subject and the number name {1} acting as the predicate. There are 46 instances (not all different) of these bidimensional expressions of unit fractions in the *Suàn Shù Shū*. They are distributed among the following patterns (the category $c_1$ includes no instances, but is known to exist in at least one other Qin-Han text):

　　　($a_1$): "Denominator + *fēn* + {1}": 24 instances presented in
　　　Sect 3–2–1
　　　($b_1$): "Denominator + *fēn* + MW + {1}": 11 instances in Sect
　　　3–2–2.
　　　($c_1$): "Denominator + *fēn* + *zhī* + {1}": no instances (Sect
　　　3–2–3).
　　　($d_1$): "Denominator + *fēn* + MW + *zhī* + {1}": 11 instances in
　　　Sect 3–2–4.

[Wenwu 2001] – and then Peng Hao (2001), [Ōkawa et al. 2006], Hu Yitao (2006) – considered that a MW was omitted between *fēn* and the numerator in some instances of the category ($a_1$) and added the MW which was implied by the context. Yang Lingrong (2008: 17–19) already argued that such additions were unnecessary; I can only emphasize that they must be rejected in the study of the actual expressions of fractions in the corpus.

**3–2–1. "Denominator + *fēn* + {1}"**
The 24 instances are:
　　　→ *sān fēn yī* for 1/3 (1 instance, on strip 3);

---

[19] This is the calculation of a surface. It also gives the conversion 1/5 *cùn* = 1/50 *chǐ* for length units (given that 1 *cùn* = 1/10 *chǐ*); the conversion is expressed as a product.

→ *sì fēn yī* for 1/4 (3 instances, on strips 3, 4 and 8);
→ *wǔ fēn yī* for 1/5 (1 instance, on strip 5);
→ *liù fēn yī* for 1/6 (2 instances, on strips 3 and 8);
→ *bā fēn yī* for 1/8 (1 instance, on strip 5);
→ *jiǔ fēn yī* for 1/9 (2 instances, on strips 3 and 8);
→ *shí fēn yī* for 1/10 (1 instance, on strip 5);
→ *shí èr fēn yī* for 1/12 (1 instance, on strip 5);
→ *shí wǔ fēn yī* for 1/15 (1 instance, on strip 6);
→ *shí liù fēn yī* for 1/16 (2 instances, on strip 5 and 9);
→ *èr shí fēn yī* for 1/20 (2 instances, on strip 6 and 9);
→ *èr shí wǔ fēn yī* for 1/25 (2 instances, on strips 6 and 8–9);
→ *sān shí fēn yī* for 1/30 (1 instance, on strip 9);
→ *sān shí liù fēn yī* for 1/36 (1 instance, on strip 9);
→ *sì shí èr fēn yī* for 1/42 (1 instance, on strip 10);
→ *sì shí jiǔ fēn yī* for 1/49 (1 instance, on strip 9);
→ *wǔ shí fēn yī* for 1/50 (1 instance, on strip 10).

### 3–2–2. "Denominator + *fēn* + MW + {1}"

The 11 instances are given in (10)–(20):

(10)　　　　　　十分　　尺　一
in *Suàn Shù Shū*　　*shí fēn*　　*chǐ*　*yī*
strip 1　　　　　　{10} *fēn*　*chǐ*　{1}
　　　　　　　　　'1/10 *chǐ*'

(11)　　　　　　廿(二十)分　尺　一
in *Suàn Shù Shū*　　*èr shí fēn*　*chǐ*　*yī*
strip 1　　　　　　{2}{10} *fēn*　*chǐ*　{1}
　　　　　　　　　'1/20 *chǐ*'

(12)　　　　　　卅(三十)分　尺　一
in *Suàn Shù Shū*　　*sān shí fēn*　*chǐ*　*yī*
strip 2　　　　　　{3}{10} *fēn*　*chǐ*　{1}
　　　　　　　　　'1/30 *chǐ*'

(13)　　　　　　八十分　　尺　一
in *Suàn Shù Shū*　　*bā shí fēn*　*chǐ*　*yī*
strip 2　　　　　　{8}{10} *fēn*　*chǐ*　{1}
　　　　　　　　　'1/80 *chǐ*'

(14)　　　　　　八分　　尺　一
in *Suàn Shù Shū*　　*bā fēn*　*chǐ*　*yī*
strip 4　　　　　　{8} *fēn*　*chǐ*　{1}
　　　　　　　　　'1/8 *chǐ*'

(15)  卌(四十)分　　尺　一
in *Suàn Shù Shū*　*sì shí fēn　chǐ　yī*
strip 4　　　{4}{10}*fēn*　*chǐ*　{1}
　　　　　'1/40 *chǐ*'

(16)  五十分　　　尺　一
in *Suàn Shù Shū*　*wǔ shí fēn　chǐ　yī*
strip 4　　　{5}{10} *fēn*　*chǐ*　{1}
　　　　　'1/50 *chǐ*'

(17)  四分　　步　一
in *Suàn Shù Shū*　*sì fēn　bù　yī*
strips 86–87　{4} *fēn*　*bù*　{1}
　　　　　'1/4 *bù*'

　　The following instances are inserted in mixed numbers.

(18)  一　錢　五分　錢　一
in *Suàn Shù Shū*　*yī　qián　wǔ fēn　qián　yī*
strip 33　{1}　*qián*　{5} *fēn*　*qián*　{1}
　　　　　'1 *qián* 1/5 *qián*'

(19)  十　斤　十二　兩　十九　朱(銖)
in *Suàn Shù Shū*　*shí　jīn　shí èr　liǎng　shí jiǔ　zhū*
strip 79　{10}　*jīn*　{10}{2}　*liǎng*　{10}{9}　*zhū*
　　　　　'10 *jīn* 12 *liǎng* 19 *zhū*

　　　　五分　朱(銖)　一
　　　　*wǔ fēn　zhū　yī*
　　　　{5} *fēn*　*zhū*　{1}
　　　　1/5 *zhū*'

(20)  七　斗　三分　升　一
in *Suàn Shù Shū*　*qī　dǒu　sān fēn　shēng　yī*
strip 119　{7}　*dǒu*　{3} *fēn*　*shēng*　{1}
　　　　　'7 *dǒu* 1/3 *shēng*'

**3–2–3. "Denominator + *fēn* + *zhī* + {1}"**
In the whole *Suàn Shù Shū*, there are no instances of the sequence *fēn zhī* uninterrupted by a measure word before {1}. This absence does not prove the pattern to be impossible; actually there is an instance of it with the expression *sān shí fēn zhī yī* 卅(三十)分之一 {3}{10} *fēn zhī* {1} on strip 0778 of *Shù* (Xiao Can 2010: 51).

### 3–2–4. "Denominator + *fēn* + MW + *zhī* + {1}"

The 11 instances of this pattern are given in (21)–(26). All are quantification phrases starting with a noun.

(21) 粺　　　　五分　　升　　之　　一
in *Suàn Shù Shū*　　*bài*　　　*wǔ fēn*　　*shēng*　　*zhī*　　*yī*
strip 100　　　　milled millet　{5} *fēn*　*shēng*　*zhī*　{1}
　　　　　　'1/5 *shēng* of milled millet'

(22) 粺米　　　　四分　　升　　之　　一
in *Suàn Shù Shū*　　*bài mǐ*　　*sì fēn*　　*shēng*　　*zhī*　　*yī*
strips 101–102　　milled millet　{4} *fēn*　*shēng*　*zhī*　{1}
4 instances　　　'1/4 *shēng* of milled millet'

(23) 毀(穀)米　　　四分　　升　　之　　一
in *Suàn Shù Shū*　　*huǐ mǐ*　　*sì fēn*　　*shēng*　　*zhī*　　*yī*
strips 102, 104　polished millet　{4} *fēn*　*shēng*　*zhī*　{1}
2 instances　　　'1/4 *shēng* of polished millet'

(24) 毀(穀)　　　　四分　　升　　之　　一
in *Suàn Shù Shū*　　*huǐ*　　　*sì fēn*　　*shēng*　　*zhī*　　*yī*
strip 103　　　polished millet　{4} *fēn*　*shēng*　*zhī*　{1}
2 instances　　　'1/4 *shēng* of polished millet'

Beware that the Chinese name for 'polished millet' is *huǐ mǐ* in (23) and *huǐ* in (24).

(25) 米　　　六　　升　　四分　　升　　之　　一
in *Suàn Shù Shū*　　*mǐ*　　*liù*　*shēng*　*sì fēn*　*shēng*　*zhī*　*yī*
strip 121　husked millet　{6}　*shēng*　{4} *fēn*　*shēng*　*zhī*　{1}
　　　　　'6 *shēng* 1/4 *shēng* of husked millet'

(26) 從(縱)　一　　步　　六分　　步　　之　　一
in *Suàn Shù Shū*　　*zòng*　*yī*　*bù*　*liù fēn*　*bù*　*zhī*　*yī*
strip 121　length　{1}　*bù*　{6} *fēn*　*bù*　*zhī*　{1}
　　　　　'a length of 1 *bù* 1/6 *bù*'

### 4. GENERIC EXPRESSIONS OF NON-UNIT FRACTIONS

There are 97 instances (not all different) of expressions of non-unit fractions stating both a numerator and a denominator in the *Suàn Shù Shū*. These bidimensional expressions of fractions are distributed among the following patterns:

> ($a_2$): "Denominator + *fēn* + Numerator": 11 instances, in Sect. 4–1.
> ($b_2$): "Denominator + *fēn* + MW + Numerator": 43 instances, in Sect. 4–2.

    ($c_2$): "Denominator + *fēn* + *zhī* + Numerator": 7 instances, in Sect. 4–3.
    ($d_2$): "Denominator + *fēn* + MW + *zhī* + Numerator": 36 instance, in Sect. 4–4.

As above with category ($a_1$), I reject the additions of a MW by [Wenwu 2001] in some instances of category ($a_2$).

### 4–1. "Denominator + *fēn* + Numerator" for non-unit fractions

The 11 instances are: *jiǔ fēn èr* for 2/9 on strip 8, *jiǔ fēn qī* for 7/9 on strip 30, on strips 22–23 there are *wǔ fēn èr* for 2/5, *liù fēn sān* for 3/6, *sān fēn èr* for 2/3, *shí fēn bā* for 8/10, *shí èr fēn qī* for 7/12, and finally 4 instances, all of them inserted in mixed numbers, given in (27)–(30); among them the 2 instances in (28) and (29) come in predicative position after a mass noun.

(27)  十二　　丰(七十)二分　　十一
in *Suàn Shù Shū*　　*shí èr*　　*qī shí èr fēn*　　*shí yī*
strip 36　　{10}{2}　　{7}{10}{2} *fēn*　　{10}{1}
'12 11/72'

(28)  粺　　　　七　斗　五分　三
in *Suàn Shù Shū*　　*bài*　　*qī*　　*dǒu*　　*wǔ fēn*　　*sān*
strip 135　　milled millet　　{7}　　*dǒu*　　{5} *fēn*　　{3}
'7 3/5 *dǒu* of milled millet'

(29)  糲　　　　二　斗　五分　二
in *Suàn Shù Shū*　　*lì*　　*èr*　　*dǒu*　　*wǔ fēn*　　*èr*
strip 136　　husked millet　　{2}　　*dǒu*　　{5} *fēn*　　{2}
'2 2/5 *dǒu* of husked millet'

(30)  四　韋(圍)　二　寸　廿(二十)五分　十四
in *Suàn Shù Shū*　　*sì*　　*wéi*　　*èr*　　*cùn*　　*èr shí wǔ fēn*　　*shí sì*
strip 154　　{4}　　*wéi*　　{2}　　*cùn*　　{2}{10}{5} *fēn*　　{10}{4}
'4 *wéi* 2 14/25 *cùn*'

### 4–2. "Denominator + *fēn* + MW + Numerator" for non-unit fractions

The 43 instances are given in (31)–(72) along with the preceding integer when there is one. Of these 2 are identical and 11 follow a noun.

(31)  十一　步　有(又)
in *Suàn Shù Shū*　　*shí yī*　　*bù*　　*yòu*
strip 84　　{10}{1}　　*bù*　　and
'11 *bù* and

|  | 九十七分 | 步 | 十(七十)九 | 步[20] |
|---|---|---|---|---|
|  | *jiǔ shí qī fēn* | *bù* | *qī shí jiǔ* | *bù* |
|  | {9}{10}{7} *fēn* | *bù* | {7}{10}{9} | *bù* |
|  | 79/97 *bù*' | | | |

| (32) | 二 | 錢 | 六十分 | 錢 | 五十七 |
|---|---|---|---|---|---|
| in *Suàn Shù Shū* | *èr* | *qián* | *liù shí fēn* | *qián* | *wǔ shí qī* |
| strip 23 | {2} | *qián* | {6}{10} *fēn* | *qián* | {5}{10}{7} |
|  | '2 *qián* 57/60 *qián*' | | | | |

| (33) | 一 | 錢 | 卅(三十)分 | 錢 | 十七 |
|---|---|---|---|---|---|
| in *Suàn Shù Shū* | *yī* | *qián* | *sān shí fēn* | *qián* | *shí qī* |
| strips 23–24 | {1} | *qián* | {3}{10} *fēn* | *qián* | {10}{7} |
|  | '1 *qián* 17/30 *qián*' | | | | |

| (34) | 金 | 三 | 朱(銖) | 九分 | 朱(銖) | 五 |
|---|---|---|---|---|---|---|
| in *Suàn Shù Shū* | *jīn* | *sān* | *zhū* | *jiǔ fēn* | *zhū* | *wǔ* |
| strip 28 | gold | {3} | *zhū* | {9} *fēn* | *zhū* | {5} |
|  | '3 *zhū* 5/9 *zhū* of gold' | | | | | |

| (35) | 七分 | 朱(銖) | 六 |
|---|---|---|---|
| in *Suàn Shù Shū* | *qī fēn* | *zhū* | *liù* |
| strip 28 | {7} *fēn* | *zhū* | {6} |
|  | '6/7 *zhū*' | | |

| (36) | 金 | 二 | 朱(銖) |
|---|---|---|---|
| in *Suàn Shù Shū* | *jīn* | *èr* | *zhū* |
| strip 28 | gold | {2} | *zhū* |
|  | '2 *zhū* | | |

|  | 六十三分 | 朱(銖) | 卌(四十)四 |
|---|---|---|---|
|  | *liùshí sān fēn* | *zhū* | *sì shí sì* |
|  | {6}{10}{3} *fēn* | *zhū* | {4}{10}{4} |
|  | 44/63 *zhū* of gold' | | |

| (37) | 六十三分 | 朱(銖) | 廿(二十)二 |
|---|---|---|---|
| in *Suàn Shù Shū* | *liù shí sān fēn* | *zhū* | *èr shí èr* |
| strip 30 | {6}{10}{3} *fēn* | *zhū* | {2}{10}{2} |
|  | '22/63 *zhū*' | | |

---

[20] The third occurrence of *bù* is superfluous and likely a copyist's mistake; this changes nothing about the classification of this fraction.

(38)    五分  錢  四
in *Suàn Shù Shū* *wǔ fēn* *qiàn* *sì*
strip 33    {5} *fēn* *qiàn* {4}
      '4/5 *qiàn*'

(39)    一 寸 六十二分  寸 卅(三十)八
in *Suàn Shù Shū* *yī* *cùn* *liù shí èr fēn* *cùn* *sān shí bā*
strip 40    {1} *cùn* {6}{10}{2} *fēn* *cùn* {3}{10}{8}
      '1 *cùn* 38/62 *cùn*'

(40)    三 寸 六十二分  寸 十四
in *Suàn Shù Shū* *sān* *cùn* *liù shí èr fēn* *cùn* *shí sì*
strips 40–41   {3} *cùn* {6}{10}{2} *fēn* *cùn* {10}{4}
      '3 *cùn* 14/62 *cùn*'

(41)    六 寸 六十二分  寸 廿(二十)八
in *Suàn Shù Shū* *liù* *cùn* *liù shí èr fēn* *cùn* *èr shí bā*
strip 41    {6} *cùn* {6}{10}{2} *fēn* *cùn* {2}{10}{8}
      '6 *cùn* 28/62 *cùn*'

(42)    ∅尺 二 寸 六十二分  寸 五十六
in *Suàn Shù Shū* ∅ *chǐ* *èr* *cùn* *liù shí èr fēn* *cùn* *wǔ shí liù*
strip 41    ∅ *chǐ* {2} *cùn* {6}{10}{2} *fēn* *cùn* {5}{10}{6}
      '[one] *chǐ* 2 *cùn* 56/62 *cùn*'

(43)    ∅尺 五 寸 六十二分  寸 五十
in *Suàn Shù Shū* ∅ *chǐ* *wǔ* *cùn* *liù shí èr fēn* *cùn* *wǔ shí*
strip 41    ∅ *chǐ* {5} *cùn* {6}{10}{2} *fēn* *cùn* {5}{10}
      '[one] *chǐ* 5 *cùn* 50/62 *cùn*'

(44)    二 斗 三 升 十一 分 升 八
in *Suàn Shù Shū* *èr* *dǒu* *sān* *shēng* *shí yī* *fēn* *shēng* *bā*
strip 48    {2} *dǒu* {3} *shēng* {10}{1} *fēn* *shēng* {8}
      '2 *dǒu* 3 *shēng* 8/11 *shēng*'

(45)    一 兩 十 朱(銖)
in *Suàn Shù Shū* *yī* *liǎng* *shí* *zhū*
strip 50    {1} *liǎng* {10} *zhū*
      '1 *liǎng* 10 *zhū*

      百卌(四十)四分  朱(銖) 九十二
      *bǎi sì shí sì fēn* *zhū* *jiǔ shí èr*
      {$10^2$}{4}{10}{4} *fēn* *zhū* {9}{10}{2}
      92/144 *zhū*'

(46)     一    錢    百一十四分     錢    ⟨七十⟩一
in *Suàn Shù Shū*    *yī*    *qián*    *bǎi yī shí sì fēn*    *qián*    *qī shí yī*
strip 57    {1}    *qián*    {10²}{1}{10}{4} *fēn*    *qián*    {7}{10}{1}
'1 *qián* 71/114 *qián*'

(47)     卌⟨四十⟩分    斗    五
in *Suàn Shù Shū*    *sì shí fēn*    *dǒu*    *wǔ*
strip 59    {4}{10} *fēn*    *dǒu*    {5}
'5/40 *dǒu*'

(48)     四    錢    八分    錢    三
in *Suàn Shù Shū*    *sì*    *qián*    *bā fēn*    *qián*    *sān*
strip 59    {4}    *qián*    {8} *fēn*    *qián*    {3}
'4 *qián* 3/8 *qián*'

(49)     八    寸    十一分    寸    二
in *Suàn Shù Shū*    *bā*    *cùn*    *shí yī fēn*    *cùn*    *èr*
strips 61–62    {8}    *cùn*    {10}{1} *fēn*    *cùn*    {2}
'8 *cùn* 2/11 *cùn*'

(50)     十八    錢    十一分    錢    九
in *Suàn Shù Shū*    *shí bā*    *qián*    *shí yī fēn*    *qián*    *jiǔ*
strip 62    {10}{8}    *qián*    {10}{1} *fēn*    *qián*    {9}
'18 *qián* 9/11 *qián*'

(51)     廿⟨二十⟩五分    錢    廿⟨二十⟩四
in *Suàn Shù Shū*    *èr shí wǔ fēn*    *qián*    *èr shí sì*
strip 64    {2}{10}{5} *fēn*    *qián*    {2}{10}{4}
'24/25 *qián*'

(52)     桼⟨漆⟩    卅⟨三十⟩七分    升    卅⟨三十⟩
in *Suàn Shù Shū*    *qī*    *sān shí qī fēn*    *shēng*    *sān shí*
strip 66    lacquer    {3}{10}{7} *fēn*    *shēng*    {3}{10}
'30/37 *shēng* of lacquer'

(53)     水    二    升    卅⟨三十⟩七分    升    七
in *Suàn Shù Shū*    *shuǐ*    *èr*    *shēng*    *sān shí qī fēn*    *shēng*    *qī*
strips 66–67    water    {2}    *shēng*    {3}{10}{7} *fēn*    *shēng*    {7}
'2 *shēng* 7/37 *shēng* of water'

(54)     七    步    卅⟨三十⟩七分    步    廿⟨二十⟩三
in *Suàn Shù Shū*    *qī*    *bù*    *sān shí qī fēn*    *bù*    *èr shí sān*
strip 68    {7}    *bù*    {3}{10}{7} *fēn*    *bù*    {2}{10}{3}
'7 *bù* 23/37 *bù*'

The nouns *jiǎn* 簡 [*bamboo strip*] and *suàn* 筭(算) [*string of coins*] occur after *fēn* in (55)–(57), and also after the name of the integer in (55) and (57). They behave in the same manner as measure words.

(55)  
in *Suàn Shù Shū*  
strip 70

二百五　　簡　　八分　　簡　　七  
*èr bǎi wǔ　jiǎn　bā fēn　jiǎn　qī*  
{2}{10²}{5}　strip　{8} *fēn*　strip　{7}  
'205 strips and 7/8 of a strip'

(56)  
in *Suàn Shù Shū*  
strip 71

八分　　簡　　一  
*bā fēn　jiǎn　yī*  
{8} *fēn*　strip　{1}  
'1/8 of 1 strip'

(57)  
in *Suàn Shù Shū*  
strip 73

十七　　筭(算)  
*shí qī　suàn*  
{10}{7}　string of coins  
'17 strings of coins

二百六十九分　　　筭(算)　　　十一  
*èr bǎi liù shí jiǔ fēn　suàn　shí yī*  
{2}{10²}{6}{10}{9} *fēn*　string of coins　{10}{1}  
11/269 of a string of coins'

(58)  
in *Suàn Shù Shū*  
strip 76

百三　　錢　　四百卅(三十)分　　錢　　九十  
*bǎi sān　qián　sì bǎi sān shí fēn　qián　jiǔ shí*  
{10²}{3}　*qián*　{4}{10²}{3}{10} *fēn*　*qián*　{9}{10}  
'103 *qián* 90/430 *qián*'

(59)  
in *Suàn Shù Shū*  
strips 80–81

水　　三　　斗[21]　　四分　　升　　三  
*shuǐ　sān　dǒu　sì fēn　shēng　sān*  
water　{3}　*dǒu*　{4} *fēn*　*shēng*　{3}  
'3 *dǒu* 3/4 *shēng* of water'

(60)  
in *Suàn Shù Shū*  
strip 84

十一　　步　　有(又)  
*shí yī　bù　yòu*  
{10}{1}　*bù*　and  
'11 *bù*

---

[21] Peng Hao (2001: 76 note 10), Hu Yitao (2006: 41) and [Ōkawa et al. 2006: 93] consider, for the sake of computational coherence, that *dǒu* 斗 as found in the text is actually a copyist's error for *shēng* 升 assuming the text means '3 and 3/4 *shēng* of water'; this correction changes nothing for our classification.

|  | 九十七分 | 步 | ꞉(七十)九 |
|---|---|---|---|
|  | *jiǔ shí qī fēn* | *bù* | *qī shí jiǔ* |
|  | {9}{10}{7} *fēn* | *bù* | {7}{10}{9} |
|  | '79/97 *bù*' | | |

(61)  　　十一分　　步　五
in *Suàn Shù Shū*　*shí yī fēn*　*bù*　*wǔ*
strip 85　　{10}{1} *fēn*　*bù*　{5}
　　　　'5/11 *bù*'

(62)  　　十一　　步　有(又)
in *Suàn Shù Shū*　*shí yī*　*bù*　*yòu*
strip 92　　{10}{1}　*bù*　and
　　　　'11 *bù*

　　　　九十八分　　步　卅(四十)七
　　　　*jiǔ shí bā fēn*　*bù*　*sì shí qī*
　　　　{9}{10}{8} *fēn*　*bù*　{4}{10}{7}
　　　　'47/98 *bù*'

(63)  　　九　步　五分　步　三
in *Suàn Shù Shū*　*jiǔ*　*bù*　*wǔ fēn*　*bù*　*sān*
strip 96　　{9}　*bù*　{5} *fēn*　*bù*　{3}
　　　　'9 *bù* 3/5 *bù*'

(64)  　　粺　　卅(三十)二分　升　九
in *Suàn Shù Shū*　*bài*　*sān shí èr fēn*　*shēng*　*jiǔ*
strip 103　　milled millet　{3}{10}{2} *fēn*　*shēng*　{9}
　　　　'9/32 *shēng* of milled millet'

(65)  　　米　　七分　升　六
in *Suàn Shù Shū*　*mǐ*　*qī fēn*　*shēng*　*liù*
strip 113　　husked millet　{7} *fēn*　*shēng*　{6}
　　　　'6/7 *shēng* of husked millet'

(66)  　　粟　　七分　升　六
in *Suàn Shù Shū*　*sù*　*qī fēn*　*shēng*　*liù*
strip 115　　unhusked millet　{7} *fēn*　*shēng*　{6}
　　　　'6/7 *shēng* of unhusked millet'

(67)  　　米　　一　升　七分　升　三
in *Suàn Shù Shū*　*mǐ*　*yī*　*shēng*　*qī fēn*　*shēng*　*sān*
strip 115　　husked millet　{1}　*shēng*　{7} *fēn*　*shēng*　{3}
　　　　'1 *shēng* 3/7 *shēng* of husked millet'

(68)   二千五十五    尺
in *Suàn Shù Shū*   *èr qiān wǔ shí wǔ*   *chǐ*
strip 149   $\{2\}\{10^3\}\{5\}\{10\}\{5\}$   *chǐ*
   '2055 *chǐ*

   卅(三十)六分   尺   廿(二十)
   *sān shí liù fēn*   *chǐ*   *èr shí*
   $\{3\}\{10\}\{6\}$ *fēn*   *chǐ*   $\{2\}\{10\}$
   20/36 *chǐ*'

(69)   ⌀韋(圍)   二   寸   廿(二十)五分   寸   十四
in *Suàn Shù Shū*   ⌀ *wéi*   *èr*   *cùn*   *èr shí wǔ fēn*   *cùn*   *shí sì*
strip 153   ⌀ *wéi*   $\{2\}$   *cùn*   $\{2\}\{10\}\{5\}$ *fēn*   *cùn*   $\{10\}\{4\}$
   '[one] *wéi* 2 14/25 *cùn*'[22]

(70)   七   寸   五分   寸   三
in *Suàn Shù Shū*   *qī*   *cùn*   *wǔ fēn*   *cùn*   *sān*
strips 153, 154   $\{7\}$   *cùn*   $\{5\}$ *fēn*   *cùn*   $\{3\}$
2 instances   '7 *cùn* 3/5 *cùn*'

(71)   從(縱)   九十七   步   有(又)
in *Suàn Shù Shū*   *zòng*   *jiǔ shí qī*   *bù*   *yòu*
strips 171–173[23]   length   $\{9\}\{10\}\{7\}$   *bù*   and
   'a length of 97 *bù* and 141/147 *bù*'

   百卅(四十)七分   步   [?][24]百卅(四十)一
   *bǎi sì shí qī fēn*   *bù*   [?] *bǎi sì shí*
   $\{10^2\}\{4\}\{10\}\{7\}$ *fēn*   *bù*   [?] $\{10^2\}\{4\}\{10\}\{1\}$

---

[22] The term *wéi* is a unit of length used for circumferences; it also appeared in (30). The shaded ⌀ signals the absence of the number name {1} before the measure word *wéi*, this also occurs in (9), (42) and (43) before *chǐ*.

[23] Strip 171 is followed by strip 173. [Wenwu 2001] and Peng Hao (2001: 117) initially had the sequence 171-172-173, but Peng Hao had changed this to 171-173-172 in the last release of the digital corpus he sent me in 2010. Hu Yitao (2006: 80, 82 note 14) and [Ōkawa et al. 2006: 1, 7 note 12] all have 171-173.

[24] I write a question mark [?] where the original shows an unclear written mark which looks like *wǔ* 五 {5}. But this would give a numerator equal to 541, which seems unlikely since it would yield an improper fraction. For the sake of coherence, the numerator should be 141. Peng Hao (2001: 121 note 39) says the digit *wǔ* 五 {5} should be corrected to *yī* 一 {1}, while Hu Yitao (2006: 82 note 14) and [Ōkawa et al. 2006: 93] say it is superfluous. This later formulation is in my view more acceptable since {1} was never used before the highest pivot of a number name in the *Suàn Shù Shū*. In any case this instance is undoubtedly an example of the pattern "Denominator's name *fēn* + MW + Numerator's name" whatever the value of the numerator.

(72)　　　　　方　　十五　　步　　卅(三十)一分　　步　　十五
in *Suàn Shù Shū*　　*fāng*　*shí wǔ*　*bù*　*sān shí yī fēn*　*bù*　*shí wǔ*
strip 185　　　　side　{10}{5}　*bù*　{3}{10}{1} *fēn*　*bù*　{10}{5}
　　　　　　　'a side of 15 *bù* 15/31 *bù*' (the length of the side)

### 4–3. "Denominator + *fēn* + *zhī* + Numerator" for non-unit fractions

The 7 instances are given in (73)–(79). Among them, 4 follow a noun and 1 follows a "Verb + OBJ" sequence. The numerical expressions in (73) and (75) are 2 of the only 4 instances of expressions of fractions with *zhī* which are not in predicative position after a noun or a phrase, the other 2 instances being those in (81) and (82) of the next section.

(73)　　　　　二千一十六分　　　　之　　百六十二
in *Suàn Shù Shū*　*èr qiān yī shí liù fēn*　*zhī*　*bǎi liù shí èr*
strip 20　　　　{2}{10³}{1}{10}{6} *fēn*　*zhī*　{10²}{6}{10}{2}
　　　　　　　'162/2016'

In (74), the expression for the fraction occurs after the pronoun object *zhī* of the verb *yuē*. The fraction is in the position of predicate with regard to the phrase *yuē zhī*; it expresses the result of a reduction.

(74)　　　　　約　　之　　百一十二分　　　之　　九
in *Suàn Shù Shū*　*yuē*　*zhī*　*bǎi yī shí èr fēn*　*zhī*　*jiǔ*
strip 20　　　　reduce　3OBJ　{10²}{1}{10}{2} *fēn*　*zhī*　{9}
　　　　　　　'reduce it [referring to 162/2016], [it is] 9/112'

In (75) the fraction is the object of a verb.

(75)　　　　　各　　受　　卅(三十)分　　之　　廿(二十)三
in *Suàn Shù Shū*　*gè*　*shòu*　*sān shí fēn*　*zhī*　*èr shí sān*
strip 26　　　　each　get　{3}{10} *fēn*　*zhī*　{2}{10}{3}
　　　　　　　'each gets 23/30'

(76)　　　　　十三　　盧唐　　四分　　之　　三
in *Suàn Shù Shū*　*shí sān*　*lútáng*　*sì fēn*　*zhī*　*sān*
strips 129–130　{10}{3}　bamboo tube　{4} *fēn*　*zhī*　{3}
　　　　　　　'13 bamboo tubes 3/4'

(77)　　　　　盾(腯)　　九分　　之　　五
in *Suàn Shù Shū*　*tú*　*jiǔ fēn*　*zhī*　*wǔ*
strip 82　　　　lard　{9} *fēn*　*zhī*　{5}
　　　　　　　'5/9 of lard'[25]

---

[25] The weigh unit *jīn* which can be deduced from the context is omitted after *fēn*. For our classification we only need to acknowledge this absence.

| (78) | | 田 | 七[26]分 | | 之 | 四 |
|---|---|---|---|---|---|---|
| in *Suàn Shù Shū* | | *tián* | *qī fēn* | | *zhī* | *sì* |
| strip 162 | | field | {7} *fēn* | | *zhī* | {4} |
| | | 'a field of 4/7'[27] | | | | |

| (79) | | 從(縱) | 廿(二十)一分 | 之 | 十六 |
|---|---|---|---|---|---|
| in *Suàn Shù Shū* | | *zòng* | *èr shí yī fēn* | *zhī* | *shí liù* |
| strip 162 | | length | {2}{10}{1} *fēn* | *zhī* | {10}{6} |
| | | 'a length is 16/21'[28] | | | |

**4–4. "Denominator + *fēn* + MW + *zhī* + Numerator" for non-unit fractions**

The 36 instances are given in (80)–(115) together with the preceding integer when the fraction is inserted in a mixed number. Note that in (81) and (82), the numerical expressions are objects of a verb, which was already the case with (75) in the previous section. The status of (115) in this matter is unclear because some characters are illegible. In all of the other 33 instances the numerical expression is inserted in a predicative clause, and the subject is a mass noun.

| (80) | | 金 | 七分 | 朱(銖) | 之 | 三 |
|---|---|---|---|---|---|---|
| in *Suàn Shù Shū* | | *jīn* | *qī fēn* | *zhū* | *zhī* | *sān* |
| strip 30 | | gold | {7} *fēn* | *zhū* | *zhī* | {3} |
| | | '3/7 *zhū* of gold' | | | | |

| (81)[29] | | 長者 | 受 | 十六 | 尺 | 有(又) |
|---|---|---|---|---|---|---|
| in *Suàn Shù Shū* | | *zhǎngzhě* | *shòu* | *shí liù* | *chǐ* | *yòu* |
| strip 55 | | elder | get | {10}{6} | *chǐ* | and |
| | | 'the elder gets 16 *chǐ* and | | | | |

| | 十八分 | 尺 | 之 | 十二 |
|---|---|---|---|---|
| | *shí bā fēn* | *chǐ* | *zhī* | *shí èr* |
| | {10}{8} *fēn* | *chǐ* | *zhī* | {10}{2} |
| | 12/18 *chǐ* ' | | | |

---

[26] The character on the strip looks more like *yī* 一 {1} than *qī* 七 {7}; but linguistic and conceptual coherence require us to chose *qī* 七 {7} to have a proper fraction. Peng Hao (2001: 115 note 5), Hu Yitao (2006: 79) and [Ōkawa et al. 2006: 20 note 4] all discussed this point. Anyway, it is of no consequence for our classification.

[27] This takes place in the calculation of one side of a rectangle when the other side and the surface are known. The textual coherence implies that a surface in square *bù* should be understood, even though the measure word is not stated. We only need to acknowledge this absence for our classification.

[28] As for the instance in (78), we only need to acknowledge the absence of a measure word to classify the expression.

[29] The fraction in (81) was already stated in (7) withdrawn from its insertion context.

21| | | | | | |
|---|---|---|---|---|---|
| (82) | 少者 | 受 | 八 | 尺 | 有(又) |
| in *Suàn Shù Shū* | *shàozhě* | *shòu* | *bā* | *chǐ* | *yòu* |
| strip 55 | younger | get | {8} | *chǐ* | and |

'the younger gets 8 *chǐ* and

| | | | | |
|---|---|---|---|---|
| | 十八分 | 尺 | 之 | 六 |
| | *shí bā fēn* | *chǐ* | *zhī* | *liù* |
| | {10}{8} *fēn* | *chǐ* | *zhī* | {6} |

6/18 *chǐ* '

| | | | | | |
|---|---|---|---|---|---|
| (83) | 粺 | 十分 | 升 | 之 | 三 |
| in *Suàn Shù Shū* | *bài* | *shí fēn* | *shēng* | *zhī* | *sān* |
| strip 98 | milled millet | {10} *fēn* | *shēng* | *zhī* | {3} |

'3/10 *shēng* of milled millet'

| | | | | | |
|---|---|---|---|---|---|
| (84) | 米 | 十五分 | 升 | 之 | 四 |
| in *Suàn Shù Shū* | *mǐ* | *shí wǔ fēn* | *shēng* | *zhī* | *sì* |
| strip 98 | husked millet | {10}{5} *fēn* | *shēng* | *zhī* | {4} |

'4/15 *shēng* of husked millet'

| | | | | | |
|---|---|---|---|---|---|
| (85) | 粟 | 廿(二十)七分 | 升 | 之 | 十 |
| in *Suàn Shù Shū* | *sù* | *èr shí qī fēn* | *shēng* | *zhī* | *shí* |
| strip 99 | unhusked millet | {2}{10}{7} *fēn* | *shēng* | *zhī* | {10} |

'10/27 *shēng* of unhusked millet'

| | | | | | |
|---|---|---|---|---|---|
| (86) | 米 | 九分 | 升 | 之 | 二 |
| in *Suàn Shù Shū* | *mǐ* | *jiǔ fēn* | *shēng* | *zhī* | *èr* |
| strip 99 | husked millet | {9} *fēn* | *shēng* | *zhī* | {2} |

'2/9 *shēng* of husked millet'

| | | | | | |
|---|---|---|---|---|---|
| (87) | 毇(穀) | 卌(四十)五分 | 升 | 之 | 八 |
| in *Suàn Shù Shū* | *huǐ* | *sì shí wǔ fēn* | *shēng* | *zhī* | *bā* |
| strip 100 | polished millet | {4}{10}{5} *fēn* | *shēng* | *zhī* | {8} |

'8/45 *shēng* of polished millet'

| | |
|---|---|
| (88) | 粟 |
| in *Suàn Shù Shū* | *sù* |
| strip 101 | unhusked millet |

'25/54 *shēng* of unhusked millet'

| | | | | |
|---|---|---|---|---|
| | 五十四分 | 升 | 之 | 廿(二十)五 |
| | *wǔ shí sì fēn* | *shēng* | *zhī* | *èr shí wǔ* |
| | {5}{10}{4} *fēn* | *shēng* | *zhī* | {2}{10}{5} |

(89)　　　　　　　　米　　　　十八分　　　升　　之　　五
in *Suàn Shù Shū*　　*mǐ*　　　　*shí bā fēn*　　*shēng*　*zhī*　*wǔ*
strip 101　　　　　husked millet　{10}{8} *fēn*　*shēng*　*zhī*　{5}
　　　　　　　　　'5/18 *shēng* of husked millet'

(90)　　　　　　　　毀(穀)米　　九分　　　升　　之　　二
in *Suàn Shù Shū*　　*huǐ mǐ*　　*jiǔ fēn*　　*shēng*　*zhī*　*èr*
strip 102　　　　　polished millet　{9} *fēn*　*shēng*　*zhī*　{2}
　　　　　　　　　'2/9 *shēng* of polished millet'

(91)　　　　　　　　麥　　　　十二分　　　升　　之　　五
in *Suàn Shù Shū*　　*mài*　　*shí èr fēn*　　*shēng*　*zhī*　*wǔ*
strip 102　　　　　wheat　{10}{2} *fēn*　*shēng*　*zhī*　{5}
　　　　　　　　　'5/12 *shēng* of wheat'

(92)　　　　　　　　米　　　　十六分　　　升　　之　　五
in *Suàn Shù Shū*　　*mǐ*　　　*shí liù fēn*　　*shēng*　*zhī*　*wǔ*
strips 102–103　　husked millet　{10}{6} *fēn*　*shēng*　*zhī*　{5}
　　　　　　　　　'5/16 *shēng* of husked millet'

(93)　　　　　　　　麥　　　　卅(三十)二分　　升　　之　　十五
in *Suàn Shù Shū*　　*mài*　　*sān shí èr fēn*　*shēng*　*zhī*　*shí wǔ*
strip 103　　　　　wheat　{3}{10}{2} *fēn*　*shēng*　*zhī*　{10}{5}
　　　　　　　　　'15/32 *shēng* of wheat'

(94)　　　　　　　　粟
in *Suàn Shù Shū*　　*sù*
strip 104　　　　　unhusked millet
　　　　　　　　　'25/48 *shēng* of unhusked millet'

　　　　　　　　　卅(四十)八分　　　升　　之　　廿(二十)五
　　　　　　　　　*sì shí bā fēn*　　*shēng*　*zhī*　*èr shí wǔ*
　　　　　　　　　{4}{10}{8} *fēn*　*shēng*　*zhī*　{2}{10}{5}

(95)　　　　　　　　粟
in *Suàn Shù Shū*　　*sù*
strip 105　　　　　unhusked millet
　　　　　　　　　'500/789 *shēng* of unhusked millet'

　　　　　　　　　七百八十九分　　　　升　　之　　五百
　　　　　　　　　*qī bǎi bā shí jiǔ fēn*　　*shēng*　*zhī*　*wǔ bǎi*
　　　　　　　　　{7}{$10^2$}{8}{10}{9} *fēn*　*shēng*　*zhī*　{5}{$10^2$}

(96)      粟      一      升
in *Suàn Shù Shū*      *sù*      *yī*      *shēng*
strips 105–106      unhusked millet      {1}      *shēng*
'1 *shēng*

     二百六十三分      之      二百卅(三十)七
     *èr bǎi liù shí sān fēn*      *zhī*      *èrbǎi sān shí qī*
     {2}{$10^2$}{6}{10}{3} *fēn*      *zhī*      {2}{$10^2$}{3}{10}{7}
237/263 *shēng* of unhusked millet'

(97)      粟      一      斗      九      升
in *Suàn Shù Shū*      *sù*      *yī*      *dǒu*      *jiǔ*      *shēng*
strip 106      unhusked millet      {1}      *dǒu*      {9}      *shēng*
'1 *dǒu* 9 *shēng*

     有(又)      二百六十三分      升      之      三
     *yòu*      *èr bǎi liù shí sān fēn*      *shēng*      *zhī*      *sān*
     and      {2}{$10^2$}{6}{10}{3} *fēn*      *shēng*      *zhī*      {3}
and 3/263 *shēng* of unhusked millet'

(98)      粟      十九      斗      有(又)
in *Suàn Shù Shū*      *sù*      *shí jiǔ*      *dǒu*      *yòu*
strip 106      unhusked millet      {10}{9}      *dǒu*      and
'19 *dǒu* and

     二百六十三分      升      之      卅(三十)
     *èr bǎi liù shí sān fēn*      *shēng*      *zhī*      *sān shí*
     {2}{$10^2$}{6}{10}{3} *fēn*      *shēng*      *zhī*      {3}{10}
30/263 *shēng* of unhusked millet'

(99)      粟
in *Suàn Shù Shū*      *sù*
strip 107      unhusked millet
'100/171 *shēng* of unhusked millet'

     百七(七十)一分      升      之      百
     *bǎi qī shí yī fēn*      *shēng*      *zhī*      *bǎi*
     {$10^2$}{7}{10}{1} *fēn*      *shēng*      *zhī*      {$10^2$}

(100)      粟      一      升      有(又)
in *Suàn Shù Shū*      *sù*      *yī*      *shēng*      *yòu*
strip 107      unhusked millet      {1}      *shēng*      and
'1 *shēng* and

        二百八十五分　　升　之　二百七(七十)五
        *èr bǎi bā shí wǔ fēn　　shēng　zhī　èr bǎi qī shí wǔ*
        {2}{10$^2$}{8}{10}{5} *fēn　shēng　zhī*　{2}{10$^2$}{7}{10}{5}
        275/285 *shēng* of unhusked millet'

(101)　　　　　粟　　　　　十七　　　升　　有(又)
in *Suàn Shù Shū*　*sù*　　　　*shí qī*　　*shēng*　*yòu*
strip 108　　　unhusked millet　{10}{7}　*shēng*　and
        '17 *shēng* and

        二百八十五分　　　升　之　百五十
        *èr bǎi bā shí wǔ fēn　　shēng　zhī　bǎi wǔ shí*
        {2}{10$^2$}{8}{10}{5} *fēn　shēng　zhī*　{10$^2$}{5}{10}
        150/285 *shēng* of unhusked millet'

(102)　　　　　粟　　　　　十七　　斗　五　升
in *Suàn Shù Shū*　*sù*　　　　*shí qī*　*dǒu*　*wǔ*　*shēng*
strip 106　　　unhusked millet　{10}{7}　*dǒu*　{5}　*shēng*
        '17 *dǒu* 5 *shēng*

     有(又)　二百八十五分　　升　之　百廿(二十)五
     *yòu*　*èr bǎi bā shí wǔ fēn　　shēng　zhī　bǎi èr shí wǔ*
     and　{2}{10$^2$}{8}{10}{5} *fēn　shēng　zhī*　{10$^2$}{2}{10}{5}
     and 125/285 *shēng* of unhusked millet'

(103)　　　　　米　　　　　卅(四十)六　石
in *Suàn Shù Shū*　*mǐ*　　　　*sì shí liù*　*shí*
strip 146　　　husked millet　{4}{10}{6}　*shí*
        '46 *shí*

        廿(二十)七分　石　之　八
        *èr shí qī fēn　shí　zhī　bā*
        {2}{10}{7} *fēn　shí　zhī*　{8}
        8/27 *shí* of husked millet'

(104)　　　　　廣　　八分　步　之　六
in *Suàn Shù Shū*　*guǎng*　*bā fēn*　*bù*　*zhī*　*liù*
strip 162　　　width　{8} *fēn*　*bù*　*zhī*　{6}
        'a width of 6/8 *bù*'

(105)　　　　　廣　　七分　步　之　三
in *Suàn Shù Shū*　*guǎng*　*qī fēn*　*bù*　*zhī*　*sān*
strip 162　　　width　{7} *fēn*　*bù*　*zhī*　{3}
        'a width of 3/7 *bù*'

| | | | | | | |
|---|---|---|---|---|---|---|
| (106) | 田 | 四分 | 步 | 之 | 二 | |
| in *Suàn Shù Shū* | *tián* | *sì fēn* | *bù* | *zhī* | *èr* | |
| strip 162 | field | {4} *fēn* | *bù* | *zhī* | {2} | |
| | 'a field of 2/4 [square] *bù*' | | | | | |

| | | | | | |
|---|---|---|---|---|---|
| (107) | 從(縱) | 百卅(三十) | | 步 | 有(又) |
| in *Suàn Shù Shū* | *zòng* | *bǎi sān shí* | | *bù* | *yòu* |
| strip 168 | length | {10²}{3}{10} | | *bù* | and |
| | 'a length of 130 *bù* and | | | | |

| | | | | |
|---|---|---|---|---|
| | 十一分 | 步 | 之 | 十 |
| | *shí yī fēn* | *bù* | *zhī* | *shí* |
| | {10}{1} *fēn* | *bù* | *zhī* | {10} |
| | 10/11 *bù*' | | | |

| | | | | | |
|---|---|---|---|---|---|
| (108) | 從(縱) | 百一十五 | | 步 | 有(又) |
| in *Suàn Shù Shū* | *zòng* | *bǎi yī shí wǔ* | | *bù* | *yòu* |
| strip 169 | length | {10²}{1}{10}{5} | | *bù* | and |
| | 'a length of 115 *bù* and | | | | |

| | | | | |
|---|---|---|---|---|
| | 廿(二十)五分 | 步 | 之 | 五 |
| | *èr shí wǔ fēn* | *bù* | *zhī* | *wǔ* |
| | {2}{10}{5} *fēn* | *bù* | *zhī* | {5} |
| | 5/25 *bù*' | | | |

| | | | | | |
|---|---|---|---|---|---|
| (109) | 從(縱) | 百五 | | 步 | 有(又) |
| in *Suàn Shù Shū* | *zòng* | *bǎi wǔ* | | *bù* | *yòu* |
| strip 170 | length | {10²}{5} | | *bù* | and |
| | 'a length of 105 *bù* and | | | | |

| | | | | |
|---|---|---|---|---|
| | 百卅(三十)七分 | 步 | 之 | 十五 |
| | *bǎi sān shí qī fēn* | *bù* | *zhī* | *shí wǔ* |
| | {10²}{3}{10}{7} *fēn* | *bù* | *zhī* | {10}{5} |
| | 15/137 *bù*' | | | |

| | | | | | |
|---|---|---|---|---|---|
| (110) | 從(縱) | 九十二 | | 步 | 有(又) |
| in *Suàn Shù Shū* | *zòng* | *jiǔ shí èr* | | *bù* | *yòu* |
| strips 172–183[30] | length | {9}{10}{2} | | *bù* | and |
| | 'a length of 92 *bù* and | | | | |

---

[30] Strip 172 is followed by strip 183, see Hu Yitao (2006: 82).

千八十九分　　　步　之　六百一十二
*qiān bā shí jiǔ fēn　　bù　zhī　liù bǎi yī shí èr*
${10^3}{8}{10}{9}$ *fēn*　*bù*　*zhī*　${6}{10^2}{1}{10}{2}$
612/1089 *bù*'

(111)　　　　　　從(縱)　八十八　　步　有(又)
in *Suàn Shù Shū*　　*zòng*　*bā shí bā*　*bù*　*yòu*
strip 175　　　　　length　${8}{10}{8}$　*bù*　and
'a length of 88 *bù* and

二千二百八十三分　　　　步　之　六百九十六
*èr qiān èr bǎi bā shí sān fēn*　*bù*　*zhī*　*liù bǎi jiǔ shí liù*
${2}{10^3}{2}{10^2}{8}{10}{3}$　*bù*　*zhī*　${6}{10^2}{9}{10}{6}$
696/2283 *bù*'

(112)　　　　　　從(縱)　八十四　　步　有(又)
in *Suàn Shù Shū*　　*zòng*　*bā shí sì*　*bù*　*yòu*
strips 177–178　　length　${8}{10}{4}$　*bù*　and
'a length of 84 *bù* and

七千一百廿(二十)九分　　　　步
*qī qiān yī bǎi èr shí jiǔ fēn*　　*bù*
${7}{10^3}{1}{10^2}{2}{10}{9}$ *fēn*　*bù*
5764/7129 *bù*'

之　五千七百六十四
*zhī*　*wǔ qiān qī bǎi liù shí sì*
*zhī*　${5}{10^3}{7}{10^2}{6}{10}{4}$

(113)　　　　　　從(縱)　八十一　　步　有(又)
in *Suàn Shù Shū*　　*zòng*　*bā shí yī*　*bù*　*yòu*
strips 180–181　　length　${8}{10}{1}$　*bù*　and
'a length of 81 *bù* and

七千三百八十一分　　　　步
*qī qiān sān bǎi bā shí yī fēn*　　*bù*
${7}{10^3}{3}{10^2}{8}{10}{1}$ *fēn*　*bù*
68??/7381 *bù*'

| | 之 | 六千八百 | | [illegible][31] |
|---|---|---|---|---|
| | *zhī* | *liù qiān bā bǎi* | | [illegible] |
| | *zhī* | {6}{10³}{8}{10²} | | [illegible] |

(114)　　　　　　廣
in *Suàn Shù Shū*　　*guǎng*
strip 183　　　　　width
　　　　　　　　'a width of

| 七 | 步 | 卅(四十)九分 | 步 | 之 | [illegible][32] |
|---|---|---|---|---|---|
| *qī* | *bù* | *sì shí jiǔ fēn* | *bù* | *zhī* | [illegible] |
| {7} | *bù* | {4}{10}{9} *fēn* | *bù* | *zhī* | [illegible] |

7 and [illegible]/9 *bù*'

(115)　　　　[illegible][33]　六十四　　　步　有(又)
in *Suàn Shù Shū*　[illegible]　*liù shí sì*　　*bù*　*yòu*
strip 183　　　[illegible]　{6}{10}{4}　　*bù*　and
　　　　　　'[illegible] 64 *bù* and

| 三百卅(四十)三分 | 步 | 之 | 二百丰(七十)三 |
|---|---|---|---|
| *sān bǎi sì shí sān fēn bù* | *bù* | *zhī* | *èrbǎi qī shí sān* |
| {3}{10²}{4}{10}{3} *fēn* | *bù* | *zhī* | {2}{10²}{7}{10}{3} |

273/343 *bù*'

### 5. LEXICALIZED EXPRESSIONS FOR 1/2, 1/3 AND 2/3

In the *Suàn Shù Shū*, the terms *bàn* 半 [half], *shǎobàn* 少半 [the smaller half] and *tàibàn* 大半 [the larger half] [34] are used as exact number names in calculations. They are lexicalized expressions of the fractions 1/2, 1/3 and 2/3; these values are revealed in (116)−(118). They can appear in expressions of mixed numbers in the order *bàn shǎobàn*, i.e. {1/2} {1/3} for 1/2+1/3 on strip 26 in (119), or in the order *shǎobàn bàn*, i.e. {1/3} {1/2} for 1/3+1/2 on strip 23. There are no other lexicalized forms for fractions in the whole *Suàn Shù Shū*[35].

The instances in (120) and (121) show that these numerals can work as verbs in "NUM + OBJ" constructions with the meaning of multiplying the value of the object by the numeral; this capability is shared by the names of integers.

---

[31] The tens and units digits are illegible and transcribed as ?? in the English translation. This does not change anything to the classification of the expression.

[32] The numerator is illegible and transcribed as [illegible] in the English translation.

[33] A noun or the upper rank digits of the integer are illegible and transcribed as [illegible] in the English translation.

[34] Only the instance of *tàibàn* on strip 8 is written 大半, the other three instances are written 泰半. I chose the reading *tàibàn* which fits the two written forms. The reading *dàbàn* would be possible for 大半 but hardly for 泰半.

[35] In Contemporary Chinese only *yī bàn* 一半 for 1/2 remains. The terms *tàibàn* (written 泰半 or 太半) or *dàbàn* 大半 are still used on occasion today but only as approximate numbers meaning *most*, no longer as exact numbers.

(116)  一半　　乘　　　一　　半　　也
in *Suàn Shù Shū*  *yī bàn*　　*chéng*　*yī*　　*bàn*　*yě*
strip 3  {1}{1/2}　multiply　{1}　{1/2}　DECL
'1/2 times 1 is 1/2

　　　　　乘　　　半　　四分　　一　　也
　　　　　*chéng*　*bàn*　*sì fēn*　*yī*　*yě*
　　　　　multiply　{1/2}　{4} *fēn*　{1}　DECL
　　　　　times 1/2 is 1/4'

(117)  少半　　　乘　　　少半
in *Suàn Shù Shū*  *shǎobàn*　*chéng*　*shǎobàn*
strip 8  {1/3}　　multiply　{1/3}
'1/3 times 1/3

　　　　　九分　　一　　也
　　　　　*jiǔ fēn*　*yī*　*yě*
　　　　　{9} *fēn*　{1}　DECL
　　　　　is 1/9'

(118)  少半　　　乘　　　大半　　　九分　　二　　也
in *Suàn Shù Shū*  *shǎobàn*　*chéng*　*tàibàn*　*jiǔ fēn*　*èr*　*yě*
strip 8  {1/3}　　multiply　{2/3}　　{9} *fēn*　{2}　DECL
'1/3 times 2/3 is 2/9'

(119)  五　　人　　分　　三　　有(又)
in *Suàn Shù Shū*  *wǔ*　*rén*　*fēn*　*sān*　*yòu*
strip 26  {5}　person　share　{3}　and
'Five people share 3 and

　　　半　　少半　　　各　　受　　　卅(三十)分　　之　　廿(二十)三
　　　*bàn*　*shǎobàn*　*gè*　*shòu*　*sān shí fēn*　*zhī*　*èr shí sān*
　　　{1/2}　{1/3}　　each　get　　{3}{10} *fēn*　*zhī*　{2}{10}{3}
　　　1/2 1/3 [a sum of three terms], each gets 23/30'
　　　　　　　　　　　　　　　　[i.e. the result of (3+1/2+1/3)÷5].

(120)  可　　半　　半　　之
in *Suàn Shù Shū*  *kě*　*bàn*　*bàn*　*zhī*
strip 17  can　{1/2}　{1/2}　3OBJ
'If it can be multiplied by 1/2 [i.e. is divisible by 2],
then multiply it by 1/2.'[36]

---

[36] This comes from a passage about the reduction of fractions on strips 17-20.

(121)  半  母     亦  半  子
in *Suàn Shù Shū*  *bàn*  *mǔ*  *yì*  *bàn*  *zǐ*
strip 19  {1/2}  denominator  also  {1/2}  numerator
'Multiply the denominator by 1/2, and
multiply the numerator by 1/2.' [37]

(122)  半  步  乘  半  步  四分一
in *Suàn Shù Shū*  *bàn*  *bù*  *chéng*  *bàn*  *bù*  *sì fēn yī*
strip 8  {1/2}  *bù*  multiply  {1/2}  *bù*  {4} *fēn* {1}
'1/2 *bù* times 1/2 *bù* is 1/4'

(123)  二  斗  泰(大)半  斗
in *Suàn Shù Shū*  *èr*  *dǒu*  *tàibàn*  *dǒu*
strip 52  {2}  *dǒu*  {2/3}  *dǒu*
'2 *dǒu* 2/3 *dǒu*'

(124)  粟     十六  斗  泰(大)半  斗
in *Suàn Shù Shū*  *sù*  *shí liù*  *dǒu*  *tàibàn*  *dǒu*
strip 88  unhusked millet  {10}{6}  *dǒu*  {2/3}  *dǒu*
'16 *dǒu* 2/3 *dǒu* of unhusked millet'

(125)  米     六  斗  泰(大)半  斗
in *Suàn Shù Shū*  *mǐ*  *liù*  *dǒu*  *tàibàn*  *dǒu*
strip 89  husked millet  {6}  *dǒu*  {2/3}  *dǒu*
'6 *dǒu* 2/3 *dǒu* of husked millet'

(126)  三分  而  乘  一
in Suàn Shù Shū  *sān fēn*  *ér*  *chéng*  *yī*
strip 3  {3} *fēn*  and then  multiply  {1}
'1/3 times 1

三分  一  也
*sān fēn*  *yī*  *yě*
{3} *fēn*  {1}  DECL
is 1/3'

(127)  七  斗  三分  升  一
in *Suàn Shù Shū*  *qī*  *dǒu*  *sān fēn*  *shēng*  *yī*
strip 119  {7}  *dǒu*  {3} *fēn*  *shēng*  {1}
'7 *dǒu* 1/3 *shēng*'

The only expressions for 1/2 in the corpus are 1 instance of {1} *bàn*, not followed by any measure word, in (116), and 46 instances of ∅ *bàn*, 12 of which are followed by a measure word, see for example (122); 33 are not, and a last instance on strip 1 is uncertain since the following characters are illegible. The regular compound "{2} *fēn*" is nowhere to be found in the text[38].

---

[37] This is from the same passage on strips 17-20.

To express 1/3 there are 24 instances of the lexicalized *shǎobàn* (15 with a measure word, 9 without), and to express 2/3 there are 4 instances of *tàibàn*, one without any measure word on strip 8, see (118); and 3 followed by the measure word *dǒu* on strips 52, 88 and 89; see (118), (123)–(125) respectively.

The regular compound *sān fēn*, i.e. {3} *fēn*, however, is found twice in the sequence {3} *fēn* {1} to express 1/3 on strips 3 and 119: see (126) without a measure word and (127) with the measure word *shēng* inserted between *fēn* and the numerator's name {1}. It is also found 17 times (already mentioned in Sect. 3–1) in monodimensional expressions of 1/3 with the numerator 1 not stated. There are three instances in expressions of 2/3: one in the sequence {3} *fēn* {2} on strip 23 (without measure word), and two on strips 138–139 with the measure word *qián* inserted between *fēn* and the numerator {2}. Therefore, in the corpus, among the 50 instances of expressions for 1/3 or 2/3, there is a choice between the lexicalized forms (28 instances) and the regular forms built with {3} *fēn* (22 instances). Any of these forms can be used to denote dimensioned quantities (weight, length, surface, etc.) and can be followed by measure words in data, calculations or results. They can also denote dimensionless coefficients in some calculations or in presentations of arithmetical procedures.

**Tab. 1: Distribution of the expressions for  
1/3 and 2/3 in the *Suàn Shù Shū***

|      | Lexicalized names *shǎobàn* and *tàibàn* | Regular forms with *sān fēn* |    |
|------|------------------------------------------|------------------------------|----|
| MW + | 18                                       | 4                            | 22 |
| MW – | 10                                       | 18                           | 28 |
|      | 28                                       | 22                           | 50 |

The distribution given in Tab. 1 shows that there is no grammatical obligation concerning the choice between lexicalized or regular items. Nevertheless, when no measure word is present, there is some preference for the unlexicalized form with {3} *fēn* since the occurrence rate of such configurations is $18/28 \times 100 \approx 64\%$. Conversely, there is preference for the lexicalized numerical items in adjectival position before a measure word since the occurrence rate of such configurations is $18/22 \times 100 \approx 82\%$; this may be because the lexicalized items yielding noun phrases are more economical than bidimensional expressions which produce predicative clauses.

## 6. CONTEXTUAL OMISSION OF THE DENOMINATOR OF A NON-UNIT FRACTION

In four passages there are series of fractions which have the same denominator, as unambiguously shown by the *context*, but this denominator is stated only in the expression of the first fraction and is understood thereafter. Abbreviations with *fēn* not preceded by the numerator's name can be construed as free reinterpretations of the item *fēn* as a noun meaning *parts* in a given partitioning.

---

[38] In Contemporary Chinese the lexicalized fraction name *yī bàn* for 1/2 can usually be replaced by the regular form {2} *fēn zhī* {1}, but not, for example, in the time expression *bā diǎn bàn* for 8:30.

In (128), taken form a passage about the taxation of pelts, the denominator 7 is stated only once. The following occurrences of *fēn* are contextually understood to designate sevenths.

(128) 犬 出 十五 錢 七分 六
in *Suàn Shù Shū* *quǎn* *chū* *shí wǔ* *qián* *qī fēn* *liù*
strips 34–35 dog exit {10}{5} qián {7} fēn {6}
'dog pelt is taxed at 15 and 6/7 *qián* [each]

貍 出 卅(三十)一 錢 ∅分 五
*lí* *chū* *sān shí yī* *qián* *∅fēn* *wǔ*
wild cat exit {3}{10}{1} qián ∅fēn {5}
wild cat pelt is taxed at 31 and 5/[7] *qián* [each]

狐 出 六十三 錢 ∅分 三
*hú* *chū* *liù shí sān* *qián* *∅fēn* *sān*
fox exit {6}{10}{3} qián ∅fēn {3}
fox pelt is taxed at 63 and 3/[7] *qián* [each]'

In (129), from another passage about the taxation of pelts, the denominator 72 is stated only once and the following occurrences of *fēn* are understood to refer to the same partitioning.

(129) 狐 出 十二 丰(七十)二分 十一
in *Suàn Shù Shū* *hú* *chū* *shí èr* *qī shí èr fēn* *shí yī*
strips 36–37 fox pay {10}{2} {7}{10}{2} fēn {10}{1}
'fox pelt is taxed at 12 11/72 [each]

貍 出 八 ∅分 卌(四十)九
*lí* *chū* *bā* *∅fēn* *sì shí jiǔ*
wild cat pay {8} ∅fēn {4}{10}{9}
wild cat pelt is taxed at 8 49/[72] [each]

犬 出 四 ∅分 十二
*quǎn* *chū* *sì* *∅fēn* *shí èr*
dog pay {4} ∅fēn {10}{2}
dog pelt is taxed at 4 12/[72] [each]'

In (130), a passage about the taxation of crops, the denominator 47 is stated only once and and the following occurrences of *fēn* are understood to refer to the same partitioning.

(130)      禾    租    四    斗    卅(四十)七分    十二
in *Suàn Shù Shū*  hé    zū    sì    dǒu    sì shí qī fēn    shí èr
strips 43–44  millet tax {4} dǒu {4}{10}{7} fēn {10}{2}
  'the tax for millet amounts to 4 *dǒu* 12/47 *dǒu*

  麥    租    三    斗    ∅分    九
  mài   zū    sān   dǒu   ∅ fēn  jiǔ
  wheat tax   {3}   dǒu   ∅ fēn  {9}
  the tax for wheat amounts to 3 *dǒu* 9/[47] *dǒu*

  荅    租    二    斗    ∅分    廿(二十)六
  dá    zū    èr    dǒu   ∅ fēn  èr shí liù
  beans tax   {2}   dǒu   ∅ fēn  {2}{10}{6}
  the tax for beans amounts to 2 *dǒu* 26/[47] *dǒu*'

In (131), the number 36 is first announced as a divisor in the calculation of a volume and then understood as the denominator of the fraction in the result.

(131)            卅(三十)六    成[39]    今
in *Suàn Shù Shū*  sān shí liù   chéng    jīn
strip 150        {3}{10}{6}    divide   now
  '36 divides, now [we get]

  二千五十五          尺    ∅分    廿(二十)
  èr qiān wǔ shí wǔ   chǐ   ∅ fēn  èr shí
  {2}{10³}{5}{10}{5}  chǐ   ∅ fēn  {2}{10}
  2055 [cubic] *chǐ* 20/[36] [cubic] *chǐ*'

## 7. CONTEXTUAL USE OF AN INTEGER NAME TO EXPRESS A DENOMINATOR

On two occasions in the corpus, the name of an integer is used to mean a fraction: see shaded ∅ in (132) and (133). The integers {7} in (132) and {4} in (133) can be understood as fractions only because the text gives the result of the calculations.

(132)            六分     乘         七∅    卅(四十)二分    一
in *Suàn Shù Shū*  liù fēn  chéng      qī ∅   sì shí èr fēn  yī
strips 9–10      {6} fēn  multiply   {7} ∅  {4}{10}{2} fēn  {1}
  '1/6 times [1/]7 is 1/42'

---

[39] The phrase "*n chéng*" (*n* 成) is an abbreviation of "*n chéng yī*" (*n* 成一) which expresses a division by *n* (Peng Hao 2001: 108 note 3), [Ōkawa et al. 2006: 29].

| (133) | 四⌀ | 乘 | 五分 | 廿(二十)分 | 一 |
|---|---|---|---|---|---|
| in *Suàn Shù Shū* | sì ⌀ | chéng | wǔ fēn | èr shí fēn | yī |
| strip 9 | {4} ⌀ | multiply | {5} *fēn* | {2}{10} *fēn* | {1} |

'[1/]4 times 1/5 is 1/20'

In the corpus, the sequence "NUM$_1$ + MW + *chéng* 乘 + NUM$_2$ (+MW)" can express the product of two lengths yielding to a surface; the unit of measurement can be dropped when it is the same for the two numbers. The instance given in (132) parallels this pattern if we consider the item *fēn* to be freely reinterpreted as a noun fitting into the measure word slot. The instance in (133) is similar except for the permutation of the two numerical expressions.

Hu Changqing (1996) cites other instances of these abbreviations in other corpora.

## 8. CASES OF ISOLATED NUMERATORS WITH SPECIAL MARKING

Two integer names preceded by *xiǎo* 小 [small] are found: *xiǎo* {5} on strip 29 and *xiǎo* {10} on strip 166. According to Peng Hao (2001: 50), the former makes sense from the context only if we interpret it as the numerator of the fraction 5/9 which is stated before in the text. As for the latter, Peng Hao (2001: 119) deduces by analogy that it must also be the numerator of a fraction which is not actually otherwise specified. Peng Hao's interpretation is quite convincing for *xiǎo* {5} but not for *xiǎo* {10} given the respective contexts.

## 9. SUMMARY OF PROMINENT FEATURES

In the *Suàn Shù Shū*, the only inseparable fraction names were on one hand the special lexicalized expressions of 1/3, 1/2 and 2/3, and on the other hand the monodimensional expressions of unit fractions built according to the pattern "Denominator's name + *fēn*" (83 instances). There were numerals which could be inserted before measure words in the same way as names for integers.

There were lexicalized forms only for 1/2, 1/3 and 2/3. Only the lexicalized form was used for 1/2 (47 instances). But to express 1/3 or 2/3, there was a choice between the lexicalized forms (28 instances) and the regular forms built with "{3} *fēn*" (22 instances). Any of these forms could be used with or without a measure word and there was no definite grammatical obligation, but a preference for the lexicalized items when a measure word was present.

Bidimensional expressions of fractions were built as predicative phrases with the unit fraction name "Denominator's Name + *fēn*" acting as subject and with the numerator's name acting as predicate. The resulting expressions were not inseparable and when a measure word was involved it was inserted right after "Denominator + *fēn*". The morpheme *zhī* was used optionally as a marker of the predicative relation. The form of the bidimensional expression of a fraction belonged to one of the four patterns defined by whether the item *zhī* was used and whether a measure word was involved. Adding the 46 instances (not all different) of bidimensional expressions of unit fractions and the 97 instances (not all different) of bidimensional expressions of non-unit fractions we get the following distribution for the total of 143 instances:

(a): "Denominator + *fēn* + Numerator": 35 instances.
(b): "Denominator + *fēn* + MW + Numerator": 54 instances.
(c): "Denominator + *fēn* + *zhī* + Numerator": 7 instances.
(d): "Denominator + *fēn* + MW + *zhī* + Numerator": 47 instances.

The item *zhī* occurred only in bidimensional expressions followed by the numerator's name and was therefore never used with the *mono*-dimensional expressions of unit fractions. Moreover the use of *zhī* was correlated with the insertion of the fraction either as the predicate in a quantification clause or as the object of a verb; Tab. 2 provides a mapping of the situation.

**Tab. 2: Bidimensional expressions of fractions in the *Suàn Shù Shū*: The item *zhī* and the insertion of fractions**

|  | *zhī* – | *zhī* + |  |
|---|---|---|---|
| Inserted – | 76 | 2 | 78 |
| Inserted + | 13 | 51 | 64 |
|  | 89 | 53 | 142 |

The characters placed before the fraction in (115) are illegible, so I do not count it here and the grand total in Tab. 2 is only 142 and not 143.

The fraction is not inserted when it occurs as the results of the calculation as in the examples (132) and (133).

The occurrence rate of *zhī* when the fraction is inserted as a predicate or an object amounts to $51/64 \times 100 \approx 80\%$. Conversely the occurrence rate of configurations without *zhī* is $76/89 \times 100 \approx 85\%$ when the fraction is not inserted. This allows us to state that the use of *zhī* inside bidimensional fraction expressions was directly correlated with the syntactical insertion of these expressions as a dependent clause used as the predicate of a quantification phrase or as the object of a verb. Readers can refer to Anicotte (2015 b) for a detailed discussion on the use of *zhī* in the expressions for fractions in Chinese.

## 10. ADDENDUM: BIDIMENSIONAL EXPRESSIONS OF PROPORTIONS

When we talk about a *fraction* of a given quantity, we assume the fraction to be one numerical item defined by a numerator and a denominator; for instance in the statement "2/3 of 9 is 6", the fraction 2/3 is one individualized number formed with the integers 2 and 3. However a numerical *proportion* between two things, or two kinds of items, can be expressed with two separate numbers.

For example, the sequence "Noun$_1$ + NUM$_1$ + Noun$_2$ + NUM$_2$" reproduced in (134) is built with two numerals both in predicative positions.

(134)  　　　　米　　　　　　一
in *Suàn Shù Shū*　*mǐ*　　　　　*yī*
strip 119　　husked millet　　{1}
　　　　　'one [part of] husked millet,

　　　　　粟　　　　　　二
　　　　　*sù*　　　　　　*èr*
　　　　　unhusked millet　　{2}
　　　　　two [parts of] unhusked millet,

　　　　　凡　　十　　斗
　　　　　*fán*　*shí*　*dǒu*
　　　　　all　　{10}　*dǒu*
　　　　　altogether 10 *dǒu*'

This expresses a proportion of one part of *husked millet* [*mǐ* 米] for two parts of *unhusked millet* [*sù* 粟] making a total volume of 10 *dǒu*. From the proportion, we can deduce that the total amount is composed of 1/3 husked millet and 2/3 unhusked millet, however these fractions 1/3 or 2/3 are not stated and their denominator 3 does not appear at all; therefore the phrase "Noun$_1$ + NUM$_1$ + Noun$_2$ + NUM$_2$" is not the expression of a fraction, but the bidimensional expressions of a proportion; therefore they should not be included in a study on fraction names.

Yang Lingrong (2008: 15–16) lists 14 examples of such expressions of proportions in the *Suàn Shù Shū* including our example (134). To these 14 examples, we can as well add this instance on strip 52:

(135)　　　　芻　　槀　　二　　石　　今
in *Suàn Shù Shū*　*chú*　*gǎo*　*èr*　*shí*　*jīn*
strip 52　　hay　　straw　{2}　*shí*　now
　　　　　'2 *shí* of hay and straw, now [there are]

　　　　　芻　　三　　而　　槀　　二
　　　　　*chú*　*sān*　*ér*　*gǎo*　*èr*
　　　　　hay　{3}　and　straw　{2}
　　　　　three [parts of] hay and two [parts of] straw'

**REFERENCES**

ANICOTTE, Rémi.
　2015 a. Chinese names for integers. In: XU Dan 徐丹 & Jingqi FU 傅京起 (eds), *Space and Quantification in Languages of China*, 117–138. Berlin: Springer.
　2015 b. Bidimensional expressions of fractions in Chinese, in *Cahiers de Linguistique Asie Orientale* (Leiden: Brill) vol.44(1): 36–56.

CHEMLA, Karine 林力娜 & Shuchun GUO 郭書春. 2004. *Les Neuf chapitres: Le classique mathématique de la Chine ancienne et ses commentaires (édition critique bilingue chinois-français)*. Paris: Dunod.

CHEMLA, Karine 林力娜 & Biao MA 馬彪. 2011. Interpreting a newly discovered mathematical document written at the beginning of the Han dynasty in China (before 157 B.C.E.) and excavated from tomb M77 at Shuihudi (睡虎地). *SCIAMVS: Sources and Commentaries in Exact Sciences* (Kyoto) vol. 12: 159–191.

CULLEN, Christopher. 2004. *The Suàn Shù Shū 筭數書 'Writings on Reckoning': A Translation of a Chinese Mathematical Collection of the Second Century BC, with Explanatory Commentary*. Cambridge UK: Needham Research Institute.

DAUBEN, Joseph W. 2008. *Suan Shu Shu*, a book on numbers and computations, English translation with commentary. *Archive for History of Exact Sciences* (Berlin: Springer) vol.62(2): 91–178.

GUO, Shuchun 郭書春. 2002. Shìlùn "*Suàn Shù Shū*" de lǐlùn gòngxiàn yǔ biānzuǎn 試論《算數書》的理論貢獻與編纂 [Preliminary discussion on the theoretical contribution and compilation of the *Suàn Shù Shū*]. *Fǎguó Hànxué* 法國漢學[Sinologie Française] (Beijing: Tsinghua University Press) vol.6: 505–537.

HU, Changqing 胡長青. 1996. Xiān Qín fēnshù biǎoshìfǎ jí qí fāzhǎn 先秦分数表示法及其发展 [Pre-Qin fraction expressions and their development]. *Gǔhànyǔ yánjiū* 古汉语研究 [Research in ancient Chinese language] (Changsha), no.32: 45–48.

HU, Yitao 胡憶濤 2006. *Zhāngjiāshān Hànjiǎn "Suàn Shù Shū" zhěnglǐ yánjiū* 張家山漢簡《算數書》整理研究 [Study on the compilation of the "*Suàn Shù Shū*" excavated in Zhāngjiāshān]. Master thesis, Southwest University (Xīnán Dàxué 西南大學), Chongqing.

JOCHI, Shigeru 城地茂. 2001. *San Sū Sho* nihon yaku《算數書》日本訳 [Japanese translation of the *Suàn Shù Shū*], in *Wasan kenkyujo kiyō* 和算研究所紀要 [Bulletin of Wasan Institute] (Tokyo), no.4: 19–46.

[Ōkawa et al. 2006]: Chōkazan Kankan *Sansūsho* Kenkyūkai 張家山《算數書》研究会 (collective). 2006. *Kan kan "San Sū Sho": Chūgoku saiko no sūgaku-sho* 漢簡《算数書》：中国最古の数学書 [The *Suàn Shù Shū*, a Han dynasty book on bamboo strips: The oldest mathematical text in China] (ed.: ŌKAWA, Toshitaka 大川俊隆). Kyōto: Hōyū shoten 朋友書店.

PENG, Hao 彭浩. 2001. *Zhāngjiāshān Hàn jiǎn* Suàn Shù Shū *zhùshì* 張家山漢簡《算術書》注釋 [Commented edition of the *Suàn Shù Shū*, a book written on bamboo strips and excavated from a Han dynasty tomb at Zhāngjiāshān]. Beijing: Kēxué Chūbǎnshè 科學出版社 [Science Press].

[Wenwu, 2000]: Jiānglíng Zhāngjiāshān Hànjiǎn zhěnglǐ xiǎozǔ 江陵张家山汉简整理小组 [Study group of the Jiangling documents] (collective). 2000. *Jiānglíng Zhāngjiāshān Hànjiǎn "Suàn Shù Shū" shìwén* 江陵张家山汉简《算数书》释文[*Transcription of bamboo "Suan Shu Shu" or A "Book of Arithmetic" from Jiangling*]. *Wénwù* 文物 [Cultural Relics] (Beijing), no.9: 78–84.

[Wenwu, 2001]: Zhāngjiāshān èr sì qī hào Hàn mù zhújiǎn zhěnglǐ xiǎozǔ 張家山二四七號漢墓竹簡整理小組 [Study group for the collation and arrangement of the bamboo strips from Zhangjiashan Han dynasty tomb no.247] (collective). 2001. *Zhāng jiā shān Hàn mù zhújiǎn (èr shí qī hào mù)* 張家山漢墓竹簡(二四七號墓) [*Bamboo strips from a Han dynasty tomb at Zhangjiashan (Tomb no.247)*]. Beijing: Wénwù Chūbǎnshè 文物出版社 [Cultural Relics Press].

XIAO, Can 肖燦. 2010. *Yuèlù Shūyuàn cáng Qín jiǎn* Shù *yánjiū* 嶽麓書院藏秦簡《數》研究 [Study of *Shù*, a Qin dynasty text written on bamboo strips and kept at the Yuèlù Academy]. PhD dissertation, Hunan University (Húnán Dàxué 湖南大學), Changsha.

XIONG Beisheng 熊北生, Wenqing YANG 杨文清, Hanqiao TAO 陶汉桥, et al. 2008. Húběi Yúnmèng Shuìhǔdì M77 fājué jiǎn bào 湖北云梦睡虎地 M77 发掘简报 [A preliminary report of the excavation at the Shuihudi tomb M77 in Yunmeng, Hubei]. *Jiānghàn kǎogǔ* 江汉考古 [*Jianghan archaeology*] (Wuhan) 2008(4): 31–37.

YANG, Lingrong 楊玲榮. 2008. *Zhāngjiāshān Hàn jiǎn shùliàng cí yǔ chēngshùfǎ yánjiū* 張家山漢簡數量詞與稱數法研究 [Study of the quantifier and the method of calling number in the Zhangjiashan bamboo slips]. PhD dissertation, East China Normal University (Huádōng Shīfàn Dàxué 華東師範大學), Shanghai.